\documentclass[12pt]{article}

\usepackage{amsmath,amssymb,amsthm}

\usepackage{epsfig}
\usepackage{arydshln}
\newcommand{\bfm}[1]{\mbox{\boldmath ${#1}$}}
\newcommand{\nonum}{\nonumber \\}
\newcommand\eq[1] {(\ref{#1})} 
\newcommand{\beqa}{\begin{eqnarray}}
\newcommand{\eeqa}[1]{\label{#1}\end{eqnarray}}
\newcommand{\beq}{\begin{equation}}
\newcommand{\eeq}[1]{\label{#1}\end{equation}}
\newcommand{\R}{\mathbb{R}}
\newcommand{\Grad}{\nabla}
\newcommand{\Div}{\nabla \cdot}

\newcommand{\Ga}{\alpha}

\newcommand{\Ge}{\epsilon}

\newcommand{\Gc}{\chi}

\newcommand{\Gk}{\kappa}

\newcommand{\Gl}{\lambda}

\newcommand{\Gm}{\mu}

\newcommand{\Gs}{\sigma}

\newcommand{\Gy}{\psi}

\newcommand{\BGs}{\bfm\sigma}

\newcommand{\BGy}{\bfm\psi}

\newcommand{\BGG}{\bfm\Gamma}

\newcommand{\BGU}{\bfm\Upsilon}

\newcommand{\BGY}{\bfm\Psi}

\newcommand{\CE}{{\cal E}}

\newcommand{\CH}{{\cal H}}

\newcommand{\CJ}{{\cal J}}

\newcommand{\CS}{{\cal S}}

\newcommand{\CU}{{\cal U}}

\newcommand{\BCE}{{\bfm{\cal E}}}

\newcommand{\BCH}{{\bfm{\cal H}}}

\newcommand{\BCQ}{{\bfm{\cal Q}}}

\newcommand{\bpm}{\begin{pmatrix}}
\newcommand{\epm}{\end{pmatrix}}

\usepackage{url}
\usepackage{hyperref}

\def\Ba{{\bf a}}
\def\Bb{{\bf b}}

\def\Be{{\bf e}}

\def\Bg{{\bf g}}
\def\Bh{{\bf h}}

\def\Bj{{\bf j}}
\def\Bk{{\bf k}}

\def\Bs{{\bf s}}

\def\Bw{{\bf w}}
\def\Bx{{\bf x}}
\def\By{{\bf y}}

\def\BA{{\bf A}}
\def\BB{{\bf B}}
\def\BC{{\bf C}}
\def\BD{{\bf D}}
\def\BE{{\bf E}}
\def\BF{{\bf F}}
\def\BG{{\bf G}}
\def\BH{{\bf H}}
\def\BI{{\bf I}}
\def\BJ{{\bf J}}
\def\BK{{\bf K}}
\def\BL{{\bf L}}
\def\BM{{\bf M}}
 
\def\BO{{\bf O}}
\def\BP{{\bf P}}
\def\BQ{{\bf Q}}
\def\BR{{\bf R}}
\def\BS{{\bf S}}
\def\BT{{\bf T}}

\usepackage{graphics}
\usepackage{amssymb}

\usepackage{color}

\title{Using subspace substitution to obtain rapidly convergent series expansions for a class of resolvents}
\author{Graeme W. Milton\thanks{Department of Mathematics, University of Utah, e-mail: milton@math.utah.edu}
}

\begin{document}
\maketitle
\vspace{2ex}
%******************************
\begin{abstract}
Following advances in the abstract  theory of composites, we develop rapidly converging series expansions about $z=\infty$ for the resolvent ${\bf R}(z)=[z{\bf  I}-{\bf P}^\dagger{\bf  Q}{\bf  P}]^{-1}$
where ${\bf  Q}$ is an orthogonal projection and ${\bf  P}$ is such that ${\bf P}{\bf  P}^\dagger$ is an
orthogonal projection. It is assumed that the spectrum of ${\bf P}^\dagger{\bf  Q}{\bf  P}$ lies within the interval $[z^-,z^+]$
for some known $z^+\leq 1$ and $z^-\geq 0$ and that the
actions of the projections ${\bf  Q}$ and ${\bf  P}{\bf  P}^\dagger$ are easy to compute. 
The series converges in the entire $z$-plane excluding the cut $[z^-,z^+]$.  It is obtained using subspace substitution, where the desired resolvent is tied to a resolvent in a larger space and ${\bf  Q}$ gets replaced by a projection $\underline{{\bf Q}}$ that is no longer orthogonal. When $z$ is real the rate of convergence of the new method matches that of the conjugate gradient method.

 %******************************
\end{abstract}

{\bf Keywords:} Resolvent, Series Expansions, Matrix Inverses, Subspace Substitution, Iterative Methods. 
\vskip 5mm
{\bf MSC codes}: 65B10, 65F10, 65Z05, 47A10, 35B27

%******************************

\section{Introduction}
Here we seek rapidly converging series expansions about $z=\infty$ of
resolvents of the form
\beq \BR(z)=[z\BI-\BP^\dagger\BQ\BP]^{-1}, \eeq{1.1}
where  
\begin{itemize}
\item $\BQ:\BCH\to\BCQ$ is an orthogonal projection 
from a Hilbert or vector space $\BCH$ onto 
$\BCQ\subset\BCH$,
\item For some subspace $\BCE\subset\BCH$, 
$\BP$ (with adjoint $\BP^\dagger$) is a mapping from  $\BCE$ to
$\BCH$ (and vice-versa),
\item $\BP$ is such that $\BGG=\BP\BP^\dagger:\BCH\to\BCH$ is an orthogonal projection onto $\BCE$, 
\item $\BP$ is such that  $\BP^\dagger\BP=\BI:\BCE\to\BCE$ is the identity on $\BCE$. 
\end{itemize}
To begin, we keep things general and consider the resolvent
\beq \BR(z)=[z\BI-\BP^\dagger\BA\BP]^{-1}, \eeq{1.1aa}
where $\BA:\BCH\to\BCH$ is Hermitian satisfying the bounds $a^-\BI\leq \BA \leq a^+\BI$ for some known real $a^+$ and $a^-$ that could be negative. 
Such resolvents appear in the abstract theory of composites and their computation using series expansions has been the subject of much attention: see \cite{Schneider:2019:PBS} and references therein. It is not the purpose of this paper to review the field. Rather the paper aims to bring to the wider numerical analysis community an approach for developing a rapidly converging series expansion of the resolvent in \eq{1.1} that has been discovered in the theory of composites. It uses the idea of subspace substitution with non-orthogonal subspaces. 

The results are based on those  in Chapter 8 of \cite{Milton:2016:ETC} and  \cite{Milton:2020:UPLVI} which in turn are based on many works that we will cite as we proceed.

By  rescaling and shifting $\BA$ and $z$, and  rescaling $\BR$, 
we can assume without loss of generality that
\beq \BI\geq \BA\geq 0.
\eeq{1.1a}
Specifically, with the replacements 
\beq \BA\rightarrow (a^+-a^-)\BA+a^-\BI,\quad z\rightarrow (a^+-a^-)z+a^-,\quad \BR\rightarrow  \BR/(a^+-a^-),
\eeq{1.aaa}
\eq{1.1} is still satisfied with $\BA$ now satisfying \eq{1.1a}.
We make the further assumption that the resolvent $\BR_\BA(a)=(a\BI-\BA)^{-1}$ can easily be computed when $a\notin [0,1]$.
In addition, we suppose that constants $z^+\leq 1$ and $z^-\geq 0$ are known such that
\beq z^+\BI\geq \BP^\dagger\BA\BP \geq  z^-\BI. \eeq{1.1b}
Multiplying \eq{1.1b} on the left by $\BP$ and on the right by $\BP^\dagger$ gives
\beq z^+\BGG \geq \BGG\BA\BGG \geq  z^-\BGG. \eeq{1.1bb}

 In the  theory of composites the case where $\BA=\BQ$ is a projection corresponds to having a two-phase medium. Then there is a simple formula for the resolvent
$\BR_\BA(a)$:
\beq \BR_\BA(a)=(a\BI-\BQ)^{-1}=\BQ/(a-1)+(\BI-\BQ)/a. \eeq{1.1bbb}
While any Hermitian operator $\BB$ with $\BI\geq\BB\geq 0$ has the factorization $\BB=\BP^\dagger\BQ\BP$ with
\beq \BQ=\bpm \BB & \BB^{1/2}(\BI-\BB)^{1/2}
\cr \BB^{1/2}(\BI-\BB)^{1/2} & \BI-\BB \epm, \quad \BP=\bpm \BI \cr  0 \epm,
\eeq{0.0}
the difficulty is in calculating the square roots of $\BB$ and $\BI-\BB$. Therefore we will assume that it is $\BQ$ and $\BP$ (or $\BGG$) which are given, not $\BB$.

The resolvent naturally appears in the solution of the equation
\beq [z\BI-\BP^\dagger\BA\BP]\By=\Bb.\eeq{1.2b}
Multiplying on the left by $\BP/z$ and noting that $\BP\By=\BGG\BP\By$, 
we get the equivalent equations
\beq [\BGG(\BI-\BA/z)]\BE=\Bh \text{   with }  
\BE=\BP\By=\BGG\BE,\quad \Bh=\BP\Bb/z=\BGG\Bh.
\eeq{1.3}
Of course, in practice one might only be interested in solving this for one value of $z$, which then can be absorbed into $\BA$,  i.e. making the replacement $\BA\rightarrow  \BA+z\BI $. However, 
keeping $z$ allows us to assess the rates of convergence of different expansions. 
So the equations become simply
\beq \BGG\BL\BE=\Bh, \quad \BGG\BE=\BE, \quad \BGG\Bh=\Bh,  \eeq{1.3a}
with
\beq \BL=\BI+(\Gs-1)\BA , \quad  \text{where   }\Gs=(z-1)/z. \eeq{1.3aa}
In a two phase conducting composite $\BA$ is the projection $\BQ=\chi\BI$, where $\chi(\Bx)$ is 1 in phase one and 0 in phase two. Then $\BL=\Gs\Gc\BI+(1-\Gc)\BI$ is the conductivity which takes the value $\Gs$ in phase one and 1 in phase two. The operator $\BGG$ is the projection onto gradients of potentials. We will formulate this example more precisely in the next section.

It is helpful to rewrite \eq{1.3a} as 
\beq  [c\BI+\BGG(\BL-c\BI)]\BE=\Bh, \eeq{1.7}
where  $c$ is a possibly complex constant. As we will see,
the choice of $c$  is typically based on bounds on the norm of certain operators. This might not be the best choice, particularly if the bounds are not sharp. An adaptive method of choosing $c$ has been discovered and implemented \cite{Schneider:2019:BBB}.

 In the theory of composites $c$ is associated with the so called reference material $\BL_0=c\BI$. Most of the analysis here extends to the case where $\BL_0$ is not a multiple of the identity but we take  $\BL_0=c\BI$ to simplify the analysis. 

Applying $\BP$ to both sides of \eq{1.7} and recalling that $\BP^\dagger\BP=\BI$ we see that
\beq  [z\BI-\BP^\dagger\BA\BP]^{-1}=\BP[c\BI+\BGG(\BL-c\BI)]^{-1}\BP^\dagger/z.
\eeq{1.7b}
Thus, we are left with evaluating the inverse of $[c\BI+\BGG(\BL-c\BI)]$ provided $c$ is such  that this inverse exists.

To define what we mean by a rate of convergence consider a series in powers
of $\Ga$:
\beq \mathbb{S}=\sum_{n=0}^\infty
\Ga^n\BC_n. \eeq{s.1}
Let us introduce the norm $\| \BO \|$ of an operator $\BO$ (that need not be self-adjoint) defined by 
\beq \|\BO\| =\max_{\Ba,\,|\Ba|=1}|\BO\Ba|, \eeq{1.2aa}
where $|\Bb|$ denotes the norm of  a vector $\Bb$ in our Hilbert space. Then an 
estimate of the error made in truncating the series \eq{s.1} at $n=m$ is
\beq  S_m=\| \mathbb{S}- \sum_{n=0}^{m-1}
\Ga^n\BC_n \|
=\|\sum_{n=m}^\infty
\Ga^{n}\BC_n \| \leq \sum_{n=m}^\infty
|\Ga^{n}| \|\BC_n \|.
\eeq{s.2}
It is the asymptotics of this latter series that determine the rate of convergence $\Gm$. The radius of convergence of the series is 
\beq R= \frac{1}{\limsup_{m\to\infty}\|\BC_m \|^{1/m}}
\eeq{s.4}
and the series will diverge if $|\Ga|>R$ and converge if $|\Ga|<R$.
We define
\beq \Gm=\frac{|\Ga|}{R}={|\Ga|\limsup_{m\to\infty}\|\BC_m \|^{1/m}}
\eeq{s.3}
as the rate of convergence.
Thus, if  $\Gm <1$ the series converges and will be fastest when $\Gm$ is small. The series will not converge if $\Gm>1$. 

For the various series studied here $\BC_n=\BC^n$ for some operator $\BC$
and the rate of convergence is simply 
$\Gm=|\Ga|\|\BC\|$.
 Then, of course, there is an iterative procedure to calculate. $\Be=\mathbb{S}\Bh$. Clearly 
$\Be_m=\sum_{n=0}^{m-1}
\Ga^n\BC^n \Bh$ satisfies the iterative relation
\beq \Be_{1}=\Bh,\quad\Be_{m+1}=\Ga\BC\Be_m+\Bh,
\,\,m=1,2,3,\ldots,
\eeq{a.4a}
and if the iterates converge, they converge to $\Be$. 

\section{A class of examples that includes the conductivity and  Schr\"odinger equations}

Equations of the form \eq{1.3a}
are prevalent throughout physics. In many cases one is interested in solving,
in the entire domain of $\R^3$ or in a unit cell of periodicity, the equation
\beq \BT(\Grad)^\dagger\BL(\Bx)\BT(\Grad)\BGY(\Bx)=\Bg(\Bx), \eeq{z.1}
where $\BT$ has a polynomial dependence on $\Grad$, and $\Bg$ is the source term.  That a host of equations in physics can be reduced to this for was recognized by Strang \cite{Strang:1986:IAMB}. While $\BT(\Grad)$ is generally linear in $\Grad $,
for some equations, such as beam and plate equations, $\BT(\Grad)$ is quadratic in $\Grad$.  Following \cite{Milton:2020:UPLI}
  let $\BE=\BT(\Grad)\BGY(\Bx)$ and choose $\Bh$ such that 
$\Bg=\BT(\Grad)^\dagger\Bh$, the latter being best done in Fourier space.
Assuming $\BF(\Grad)=\BT(\Grad)^\dagger\BT(\Grad)$
has an inverse, easiest to compute in Fourier space, we introduce the projection
\beq \BGG(\Grad)=\BT(\Grad)(\BF(\Grad))^{-1}
\BT(\Grad)^\dagger,
\eeq{z.2}
and then \eq{1.3a} is satisfied with 
$\Bh=\BT(\Grad)^\dagger(\BF(\Grad))^{-1}$ and we can recover $\BGY$ using 
\beq  \BGY=(\BF(\Grad))^{-1}\BT(\Grad)^\dagger\BE.
\eeq{z.3}
The operators $(\BF(\Grad))^{-1}$ and hence $\BGG(\Grad)$ are non-local in real space. Their action is best evaluated in Fourier space where $\Grad$ is replaced by $i\Bk$, where $\Bk$ is the wave vector. In the context of these equations bounds $z^+$ and $z^-$ on the spectrum of $\BA$ can be obtained by using, as appropriate, quasi-convex functions \cite{Dacorogna:2007:DMC}, A-quasi-convex 
functions \cite{Braides:2000:AQR}, or their generalization of $Q^*$-convex operators \cite{Milton:2020:UPLV}.
%These equations can also be expressed in the form %\eq{1.1aa} with
%\beq \BP(\Grad)=\BT(\Grad)(\BF(\Grad))^{-1/2},
%\eeq{z.3a}
%satisfying  $\BP^\dagger\BP=\BI$ and 
%$\BP\BP^\dagger=\BGG$ . However, not that the 
%square root of $\BF^{-1}$ as opposed to  $\BF^{-1}$ does %not enter the equations \eq{1.3a}.

The prototypical example is the conductivity equation:
\beq \Bj'(\Bx)=\BGs(\Bx)\Be(\Bx)-\Bs(\Bx),\quad \BGG\Be=\Be,\quad\BGG\Bj'=0, \quad\text{with}\quad \BGG=\Grad(\Grad^2)^{-1}\Grad\cdot,
\eeq{pt1}
where $\BGs(\Bx)$ is the conductivity tensor, while $\Div\Bs$, $\Bj=\Bj'+\Bs$, and $\Be$ are the current source, current,
and electric field, and $(\Grad^2)^{-1}$ is the inverse Laplacian
(there is obviously considerable flexibility in the choice of $\Bs(\Bx)$, the
only constraints being square integrability and that
$\Div\Bs$ equals the current source).
An interesting twist is that we write the equations in this form,
rather than in the more conventional form involving $\Bj$ directly. This is exactly
what we need to keep the left hand side of the constitutive law divergence free.
As current is conserved, $\Div\Bj=\Div\Bs$, implying $\Div\Bj'=0$, which is clearly equivalent to
$\BGG\Bj'=0$. To show that $\Be=\Grad(\Grad^2)^{-1}\Grad\cdot\Be$, we let $V$ be the solution of Poisson's equation
$\Grad^2V=-\Div\Be$ (with $V(\Bx)\to 0$ as $\Bx\to\infty$), i.e. $V=-(\Grad^2)^{-1}\Div\Be$.
Then substitution gives $\Be=-\Grad V=\Grad(\Grad^2)^{-1}\Div\Be=\BGG\Be$.
These steps are much easier done in Fourier space, where $\BGG(\Bk)=\Bk(\Bk\cdot\Bk)^{-1}\Bk^T$. Upon
identifying $\BL$ with $\BGs$ and $\BE$ with $\Be$, the application of $\BGG$ to the first equation in \eq{pt1} gives \eq{1.3a} with $\Bh=\BGG\Bs$. For a two phase conducting composite with isotropic (scalar) conductivities
 we can rescale, without loss of generality, the conductivities
so phase 2 has conductivity 1. Then $\BL=\BI+(\Gs-1)\BQ$
where $\Gs$ is the conductivity of phase 1 and $\BQ=\Gc\BI$
where $\Gc$ is the characteristic function taking the value 1 in phase 1 and 0 in phase 2. 

Another example is the Schr\"odinger equation. With $m$ denoting the mass of the electron and $\hbar$ denoting Planck's constant $h$ divided by $2\pi$, the time independent Schr\"odinger equation for the wave function $\Gy$ of a single electron in an infinite domain with source term $g(\Bx)$, potential $V(\Bx)$ and energy $E_0$ takes the form
\beq [V(\Bx)-E_0]\Gy(x)-\frac{\hbar}{2m}\Grad^2\Gy=g(\Bx).
\eeq{z.4}
It can be re-expressed as 
\beq
\underbrace{\bpm -\Grad &  1 \epm}_{\BT(\Grad)^\dagger}\underbrace{\bpm \hbar\BI/2m & 0 \cr 0 &  V-E_0 \epm}_\BL
\underbrace{\bpm \Grad  \cr 1 \epm }_{\BT(\Grad)}\BGy
=g,
\eeq{z.5}
and in Fourier space
\beq \BGG= \frac{1}{1+|\Bk|^2}\bpm \Bk\otimes\Bk & i\Bk \cr -i\Bk & 1 \epm,
\quad \BE=\bpm i\Bk  \cr 1 \epm \BGy,
\quad\Bh=\frac{1}{1+|\Bk|^2}\bpm i\Bk  \cr 1 \epm g.
\eeq{z.6}
Note that $\BL$ will be bounded only if $V(\Bx)$ is bounded and coercive if only if $V(\Bx)-E_0\geq \Ge$ for some $\Ge>0$. 

Now suppose
$V(\Bx)$ only takes two values,
\beq V(\Bx)=V_0+\Gc(\Bx)V_1, \eeq{z.7}
where $\Gc(\Bx)$ takes either the value $1$ or zero and
$V_0-E_0$  and $V_0+V_1-E_0$ are both positive ($V_1$ could be negative).
We
rewrite \eq{z.4} as
\beq
\underbrace{\bpm -\Grad \sqrt{\hbar/2m}  &  \sqrt{V_0-E_0} \epm}_{\BT(\Grad)^\dagger}\underbrace{\bpm \BI & 0 \cr 0 &  1+\Gc V_1/(V_0-E_0) \epm}_{\BL}
\underbrace{\bpm \Grad \sqrt{\hbar/2m}  \\ \sqrt{V_0-E_0}\epm}_{\BT(\Grad)}\Gy=g, 
\eeq{z.7a}
where we have replaced the previous definitions of $\BT(\Grad)$ and $\BL$.  We now change the definitions of  $\BGG$, $\BE$ , and $\Bh$ in Fourier space to
\beqa \BGG&=& \frac{1}{V_0-E_0+|\Bk|^2\hbar/2m}\bpm (\Bk\otimes\Bk)\hbar/2m & i\Bk\sqrt{(V_0-E_0)\hbar/2m}\cr -i\Bk\sqrt{(V_0-E_0)\hbar/2m} & V_0-E_0\epm, \nonum
\BE & = &\bpm \Grad \sqrt{\hbar/2m}  \\ \sqrt{V_0-E_0}\epm \Gy,
\quad
\Bh= \frac{1}{V_0-E_0+|\Bk|^2\hbar/2m}
\bpm \Grad \sqrt{\hbar/2m}  \\ \sqrt{V_0-E_0}\epm g,
\eeqa{z.8}
where $\BGG$ projects onto the subspace of fields having the same form as  $\BE$ as $\Gy$ varies.

Then the equations still take the form \eq{1.3a} where now
\beq \BL=[\BI-\BQ/z],\text{    with   }\BQ= \bpm 0  & 0 \cr 0 &  \Gc \epm,\quad\quad
z=(E_0-V_0)/V_1.
\eeq{z.9}

A plethora of other equations that can be expressed in the form \eq{1.3a} are presented in 
\cite{Milton:2020:UPLI, Milton:2020:UPLII, Milton:2020:UPLIII, Milton:2020:UPLIV}.

\section{Some well known series expansions}
\setcounter{equation}{0}

The simplest series expansion of the resolvent \eq{1.1} is of course
\beq \BR=\sum_{n=0}^\infty(\BP^\dagger\BA\BP)^n/z^{n+1}
=\sum_{n=0}^\infty \BP^\dagger(\BA\BGG)^n\BP/z^{n+1},
\eeq{1.2}
where we have used the fact that $\BP\By=\BGG\BP\By$.
As is well known
one can estimate the error in taking just $m$ terms in the series expansion
\beqa \|\BR-\sum_{n=0}^{m-1}(\BP^\dagger\BA\BP)^n/z^{n+1}\| &=&
\|\sum_{n=m}^{\infty}(\BP^\dagger\BA\BP)^n/z^{n+1}\| \nonum
&  \leq  &\sum_{n=m}^{\infty}\|(\BP^\dagger\BA\BP)\|^n/|z|^{n+1} \nonum
& \leq & \sum_{n=m}^{\infty}|z^+/z|^n/|z|
=|z^+/z|^m/(|z|-z^+).
\eeqa{1.2a}
Thus, the rate of convergence for the series in powers of $1/z$ is at least $\Gm_0=z^+/z$.

An improvement is the Richardson method which expands the inverse on the right side of  \eq{1.7b} giving
\beqa [c\BI+\BGG(\BL-c\BI)]^{-1} &= & \sum_{n=0}^\infty [(\BGG(\BL/c-\BI)]^n/c
= \sum_{n=0}^\infty [(\BGG(\BL/c-\BI)\BGG]^{n-1}(\BL/c-\BI)/c.
\eeqa{1.8}
Such series \cite{Brown:1955:SMP, Kroner:1972:SCM} have played an important role in the theory of composites.
In that context $\BGG$ imposes the differential constraints of the fields and
acts locally in Fourier space while $\BL$ represents the material tensor field
entering the constitutive relation and it acts locally in 
real space. Particularly important was the recognition by Moulinec and 
Suquet \cite{Moulinec:1994:FNM}
 that fast Fourier transforms could be used to bounce back and forth between real space and Fourier
space when computing by iteration the action of the series in \eq{1.8} on a field. Their numerical method  triggered a multitude of further developments, see \cite{Schneider:2019:PBS} and references therein, including
the extension to non-linear equations in composite materials as reviewed in \cite{Schneider:2021:RNF}. 
Observe that 
\beq \|\BL/c -\BI\| \leq \max\{|\Gs-c|/|c|, |1-c|/|c|\}. \eeq{1.9}
The right side is minimized  when $c=(1+\Gs)/2$ , which is the choice made by Moulinec and 
Suquet \cite{Moulinec:1994:FNM}. In this case, the series expansion becomes
\beq  [c\BI+\BGG(\BL-c\BI)]^{-1}=2\sum_{n=0}^\infty s^n[\BGG(2\BA-\BI)]^n/(1+\Gs)
\text{  with } s=\frac{\Gs-1}{\Gs+1}.
\eeq{1.10}
Then, since $\|2\BA-\BA\|\leq 1$, the rate of convergence is at least $\Gm_1= |s|$.
If we only take $m$ terms in the series then the error can be bounded:
\beqa  \| [c\BI+\BGG(\BL-c\BI)]^{-1}-2\sum_{n=0}^{m-1} s^n[\BGG(2\BA-\BA)]^n/(1+\Gs)\| 
&= &\| 2\sum_{n=m}^\infty s^n[\BGG(2\BA-\BI)]^n/(1+\Gs) \| \nonum
&\leq & 2\sum_{n=m}^\infty |s|^n\left(\|\BGG\| \|2\BA-\BI\|\right)^n/|1+\Gs| \| \nonum
&\leq & 2\sum_{n=m}^\infty |s|^n/|1+\Gs| \nonum
& = & 2|s|^{m}/[(1-|s|)|1+\Gs|].
\eeqa{1.10a}
Although it will not concern us here, the rate of convergence of this series with this value of $c$ is better than that indicated by the bound \eq{1.10a}
when $z^-\neq 0$. To see this,  \eq{1.1bb} implies
\beq
|s|\|\BGG(2\BA-\BI)\BGG]\|\leq |s|\max\{|2z^- -1|,|2z^+ -1|\}
\eeq{1.10ab}
and thus the rate of convergence is bounded above by  the right hand side of 
\eq{1.10ab}.
%\beq \| 2\sum_{n=m}^\infty s^n[\BGG(2\BA-\BI)]^n/(1+\Gs) \| 
%\leq 2\sum_{n=m-1}^\infty |s|^n\|\BGG(2\BA-\BI)\BGG]\|^n|\|2\BA-\BI\|/(1+\Gs)
%\eeq{1.10b}
%and \eq{1.1bb} implies $ \|\BGG(2\BA-\BI)\BGG]\|=1-2z^-$
%giving 
%\beq   \| [c\BI-\BGG(\BL-c\BI)]^{-1}-2\sum_{n=0}^{m-1} s^n[\BGG(2\BA-\BI)]^n/(1+\Gs)\| 
%\leq 2|1+\Gs||s(1-2z^-)|^{m+1}/(1-|s(1-2z^-)|)
%\eeq{1.10c}
% So the rate of converge is  at least|s(1-2z^-)|.

If both $z^+$ and $z^-$ are known, further improvements can be obtained by adjusting $c$ to minimize the norm of  the operator $\BGG(\BL/c-\BI)\BGG$
appearing on the right side of \eq{1.8} results in a series with an even better bound on the rate of convergence. Specifically, using \eq{1.1bb} to bound this norm gives
\beq \|\BGG(\BL/c-\BI)\BGG\|=\|(1-c)\BGG-\BGG\BA\BGG/z\|/|c| 
\leq \max\{|(1-c)-z^+/z|/|c|, |(1-c)-z^-/z|/|c|\}.
\eeq{1.10d}
The maximum occurs when both expressions in the maximum are equal, i.e. for
\beq c=1-(z^+ + z^-)/(2z) ,
\eeq{1.10e}
and in this case we get
\beq  \|\BGG(\BL/c-\BI)\BGG\|\leq \frac{z^+-z^-}{|2z-z^- -z^+|}.
  \eeq{1.10f}
So the rate of convergence is at least 
\beq \Gm_2=(z^+-z^-)/|2z-z^+-z^+|.
\eeq{1.10g}
To compare this with the convergence rate bound $\Gm_1=|s|$ we first note that the first series will converge only if $|s|<1$.  Writing $\Gm_2$ in terms of $s$ gives
\beq \Gm_2=\frac{(z^+-z^-)|s|}{|s(1-z^+-z^+)-1|}.
\eeq{1.10h}
As $s$ varies in the unit disk the minimum of  $|s(1-z^+-z^+)-1|$ occurs at 
$s=1$  and so we get 
\beq \Gm_2\leq \frac{|s|(z^+-z^-)}{z^+ +z^-}
\leq \Gm_1. \eeq{1.10i}

\section{Motivating the strategy}
\setcounter{equation}{0}
Suppose we wanted to find an expansion that converges as rapidly as possible of  $f(z)$
about $z=\infty$ given that $f(z)$ has singularities distributed all along the interval $[z^-, z^+]$. The strategy is to map the region of the complex plane outside this interval to the strict interior of the unit disk. Thus we want the image of $[z^-, z^+]$ to wrap entirely around the boundary of the unit disk. First we make the fractional linear transformation from the complex $z$-plane to the complex $ \underline{\sigma}$-plane where
\beq \underline{\sigma}=\frac{z-z^-}{z-z^+}.  \eeq{m.2}
This maps  $[z^-, z^+]$  to the entire negative real $ \underline{\sigma}$ axis, and moves the expansion point from $z=\infty$ to  $ \underline{\sigma}=1$. Next we introduce the mapping
\beq w=\sqrt{\underline{\sigma}}, \eeq{m.3}
which takes the cut complex  $ \underline{\sigma}$-plane to the right half of the $w$-plane.
Now the singularities are distributed along the entire negative imaginary axis and the expansion
point remains at $w=1$.  Finally, we take
\beq  v=\frac{w-1}{w+1}, \eeq{m.4}
which maps the right half of the $w$-plane to the unit disk in the $v$-plane and moves the expansion point from $w=1$ to the origin in the $v$-plane. Putting all these transformations together we see that 
\beq
v=\left(
\frac{\sqrt{\frac{z-z^-}{z-z^+}}-1}
{\sqrt{\frac{z-z^-}{z-z^+}}+1}
\right)
\eeq{m.5}
is the desired transformation. We will seek an expansion of the resolvent in powers of $v$.

Putting together the inverse transformations,
\beq w=\frac{1-v}{v+1} , \quad \underline{\sigma}=w^2, \quad 
z=\frac{z^+\underline{\sigma}-z^-}{\underline{\sigma}-1},
\eeq{m.6}
we obtain 
\beq z=z(v)\equiv\frac{\frac{z^+(1-v)^2}{(v+1)^2}-z^-}{\left(\frac{1-v}{v+1}\right)^2-1}.
\eeq{m.7}
This may be substituted in $f(z)$ and the resulting function
$f(z(v))$ can then be expanded in powers of $v$
to obtain a rapidly converging expression for $f$ with the expansion point at $v=0$ corresponding to the expansion point of $f(z)$ at $z=\infty$. 

\section{Obtaining the desired series expansions}
\setcounter{equation}{0}
\subsection{An accelerated series expansion that does not require knowledge of the spectral bounds $z^+$ and $z^-$}
	
The initial stage in obtaining a  series expansion with a generally improved rate of convergence is to obtain one where $\sqrt{\sigma}$ naturally enters. (Later we will replace $\Gs$ by an appropriately defined $\underline{\Gs})$). 
The series expansions developed in this initial stage were introduced for the conductivity equations in composites \cite{Eyre:1999:FNS} and then extended to elasticity \cite{Michel:2001:CSL} and more general equations (Section 14.9 in \cite{Milton:2002:TOC}). Bounds on the convergence of this accelerated scheme were derived in this latter paper. 
Later it was recognized \cite{Schneider:2019:PBS} that the accelerated iteration scheme is a special case of the Douglas Rachford splitting method \cite{Douglas:1956:NSH}. As Pierre Suquet noted  (private communication, July 2024),  another special case of the Douglas Rachford splitting method is the augmented Lagrangian iteration scheme \cite{Michel:2000:CMB} previously found  to have a close connection 
\cite{Moulinec:2014:CTA} with the accelerated scheme
\cite{Eyre:1999:FNS}. Bounds on the convergence of this splitting method have been derived
\cite{Giselsson:2017:LCM}.

Our starting point is again the inverse on the left side of \eq{1.8} but now we change the sign of $c$, which is equivalent to taking a reference tensor $\BL_0$ that is negative definite.
Then, following \cite{Milton:1990:RCF},  we have
\beq [-c\BI+\BGG(\BL+c\BI)]^{-1}=(\BL+c\BI)^{-1}[-c(\BL+c\BI)^{-1}+\BGG]^{-1}
=2(\BL+c\BI)^{-1}(\BK-\BGU)^{-1},
\eeq{x.1}
where
\beq \BK=(\BL-c\BI)(\BL+c\BI)^{-1},\quad \BGU=\BI-2\BGG.
\eeq{x.2}
Since $\BGU$=$\BGU^{-1}$ we can rewrite this as
\beq [-c\BI+\BGG(\BL+c\BI)]^{-1}=-2(\BL+c\BI)^{-1}(\BI-\BGU\BK)^{-1}\BGU.
\eeq{x.3}
As we are interested in the action of this inverse on $\Bh=\BGG\Bh=-\BGU\Bh$ we can
replace  the last $\BGU$ in \eq{x.3} by simply $-\BI$. 

Our assumption that the resolvent $(a\BI-\BA)^{-1}$ can easily be computed when $a\notin [0,1]$ allows us to easily compute $\BK$ when $(1+c)/(1-\Gs)\notin [0,1]$ which is guaranteed to be the case when $c>0$ and $\Gs$ is not on the negative real axis. To bound $\|\BK\|$
we first consider the special case where $\Gs$ is real with $\Gs>1$ and $c$ is real and positive with $c\in [1,\Gs]$. Then the norm of $\BK$ as $\BA$ varies with $0\leq\BA\leq\BI$ is maximized when $\BA=0$
or $\BA=\BI$, giving
\beq \|\BK\|\leq \max\{\left |\frac{\Gs-c}{\Gs+c}\right |, \left |\frac{1-c}{1+c}\right | \},
\eeq{1.13}
and the bound is tightest when $c$ is chosen so the first expression in the brackets equals the second expression, that is
at $c=\sqrt{\Gs}$. With this value of $c$, and using the fact that $\|\BGU\|=1$  the series expansion comes from
\beq (\BI-\BGU\BK)^{-1}=\sum_{n=0}^\infty(\BGU\BK)^n\text{   where  }
\|\BGU\BK\|\leq \|\BGU\| \|\BK\|\leq |u| \text{   with  } u=\frac{\sqrt{\Gs}-1}{\sqrt{\Gs}+1},
\eeq{1.14}
giving
\beq 
[-c\BI+\BGG(\BL+c\BI)]^{-1}\BGU= -2(\BL+c\BI)^{-1}\sum_{n=0}^\infty(\BGU\BK)^n.
\eeq{1.14AAA}
So the rate of  convergence is  $\Gm_3=|u|$.
When $\Gs$ is complex then with $c$ and $k$ kept at $\sqrt{\Gs}$ and $1/2$
the bound on the norm of $\BK$
needs to be replaced by
\beq  \|\BK\|=\max_{\Gl\in[0,1]} |t(\Gl)|\text{      where     }
t(\Gl)=1-\frac{2\sqrt{\Gs}}{\sqrt{\Gs}+1+(\Gs-1)\Gl},
\eeq{1.16}
in which $\Gl$ is a possible eigenvalue of $\BA$. As $t(\Gl)$ is a fractional linear transformation of $\Gl$ it inscribes a circular arc in the complex plane joining
$t(0)$ with $t(1)=-t(0)$ as $\Gl$ is varied between 0 and 1, that when extended passes through $t=1$ at $\Gl=\infty$. The maximum in \eq{1.16} is attained 
at $\Gl=0$ and $\Gl=1$  if and only if $t=1$ lies on or outside  the circle in the complex $t$ plane
centered at $t=0$ and having radius $|t(1)|=-|t(0)|$ , i.e. when 
$|t(1)|=-|t(0)|\leq 1$. This is satisfied when $\Gs$ is in the right half of the complex plane and if this holds then the rate of convergence remains $\Gm_3=|u|$.

The rate of convergence $\Gm_3$  is better than $\Gm_1=|s|$ with $s=(\Gs-1)/(\Gs+1)$. To see this, we first express $\sqrt{\Gs}$ 
in terms of $u$: $\sqrt{\Gs}=(1+u)(1-u)$ and then substitute this 
back in the formula for $s$ giving
\beq s=\frac{2u}{u^2+1}. \eeq{1.16aa}
So if $u$ and hence $u^2$ lie in the unit disk then $|u^2+1|\leq 2$
implying that $\Gm_3\leq \Gm_1$.

While the rate of convergence $\Gm_3$ is generally better than the rate of convergence $\Gm_2$, i.e. $\Gm_3<\Gm_2$ this is not always the case. In particular, the
last series might not converge, i.e. $\Gm_3=1$, when $\Gs$ is on the negative real axis, while $\Gm_2=1$ when $\Gs$ is on a circle symmetric about the negative real axis intersecting it at the points $1-1/z^-$ and $1-1/z^+$ .  A comparison of the convergence of the 
two series \eq{1.10} and \eq{1.14A} for a model example in the theory of composites is given in \cite{Moulinec:2018:CIM}.

As shown in \cite{Milton:2002:TOC} there is an alternative approach which can be reduced to essentially the same  expansion when $\BA=\BQ$. It offers more flexibility that could result in faster convergence, 
involving  not only a reference operator 
$\BL_0$, now positive definite but not necessarily proportional to $\BI$, but also an additional operator $\BM$. Then one makes a manipulation of the equations into one involving the 
inverse $(\BI-\BGU\BK)^{-1}$ where now
\beq \BK=[\BI+(\BL-\BL_0)\BM]^{-1}(\BL-\BL_0), \quad \BGU=\BM-\BGG'\quad \text{with}
\quad\BGG'=\BGG(\BGG\BL_0\BGG)^{-1}\BGG,
\eeq{y.1}
in which the last inverse needs to be taken on the space on which $\BGG$ projects. $\BGG'$ can alternatively be defined through its action: 
$\BE'=\BGG'\BP$ if and only if $\BGG\BE=\BE$ and $\BGG(\BP-\BL_0\BE')=0$.
$\BM$ and $\BL_0$ need to be such that both $\BK$ and $\BGG'$ can be easily computed. 
This manipulation, introduced in \cite{Milton:1990:RCF}, and its generalizations  are key to developing alternative formulas for the effective tensors of laminates \cite{Milton:1990:CSP, Zhikov:1991:EHM}.
Combined with the associated series expansions, it enables the development of the general theory of exact relations for hierarchical laminates \cite{Grabovsky:1998:EREa} and for arbitrary microstructures
\cite{Grabovsky:2000:ERE} (see also Chapter 17 in \cite{Milton:2002:TOC}
and the book \cite{Grabovsky:2016:CMM}). Exact relations are identities that hold irrespective of the microstructure provided the material tensors lie on suitable manifolds.

When $\BL_0=c\BI$ (now positive definite) and $\BM=\BI/2c$ the formulas for $\BK$ and $\BGU$ , with $\BK$ divided by $2c$ and $\BM$ multiplied by $2c$ match those in
\eq{x.2}. So we obtain the same series expansion.

\subsection{The new series with a faster convergence rate when $\BA$ is a projection}

We now specialize to the case of interest where $\BA=\BQ$ is an orthogonal  projection.
Our goal is  to find a series expansion, with a rate of convergence at  most $\Gm_4$,
which is superior in the sense that $\Gm_4$ is always less than or equal to $\Gm_2$ and $\Gm_3$ (and hence $\Gm_1$ since $\Gm_2\leq \Gm_1$). Then,  the
resolvent $\BR_\BQ{\Ba}= (a\BI-\BQ)^{-1}$ is given by \eq{1.1bbb}
and hence $\BK$, given by \eq{x.2}, 
simplifies to $u(2\BQ-\BI)$. Consequently , the series expansion \eq{1.14AAA} with $c=\sqrt{\Gs}$
 becomes
 \beq [\BI-{\BGG}{\BQ}/{z}]^{-1}= \BH\sum_{n=0}^\infty u^n\BD^n\text{     with  }
  u=\frac{\sqrt{{\sigma}}-1}
 {\sqrt{{\sigma}}+1},
 \eeq{1.14A}
where
 \beq \BH
 =2(\BL+\sqrt{\Gs}\BI)^{-1}=2[\BQ/\sqrt{\Gs}+(\BI-\BQ)]/(1+\sqrt{\Gs}), \quad
 \BD=(2{\BQ}-\BI)(\BI-2\BGG), 
 \eeq{1.16A}
 and has the same rate of convergence $\Gm_3=|u|$.
 
 The next step to obtain a  series expansion with an improved rate of convergence is  to find one where a fractional linear transformation mapping  the interval $[z^-,z^+]$ in the complex $z$ plane to the negative real axis naturally enters. 
 The analysis here is based on that given in Chapter 8 of \cite{Milton:2016:ETC}
 
  To start we consider the following linear algebra problem: given $s_1,s_2,s_3$ and $t$, solve the matrix equation,
 \beq
 \bpm
 J\\0 \\ J_2
 \epm = \bpm E\\E_1 \\ 0\epm-\frac{1}{\underline{z}}\underbrace{\bpm s_1^2 & s_1s_2 & s_1s_3 \\ s_1s_2 &
 	s_2^2 & s_2s_3 \\  s_1s_3 & s_2s_3 & s_3^2 \epm}_{\BS}\bpm E\\E_1\\0\epm-\bpm h\\ 0 \\0\epm,
 \eeq{4.1}
 for $J$ in terms of $E$ and $h$. We will ultimately allow for $s_1,s_2$ and $s_3$, that are either real or purely imaginary, chosen with
 \beq s_3^2=1-s_1^2-s_2^2, \eeq{4.1a}
 to ensure that $\BS$ is a projection matrix, though not self-adjoint in our application implying that $\BS$ and $\BI-\BS$ project onto subspaces that are not orthogonal. 
 The significance of \eq{4.1} is that it corresponds to a problem in the abstract theory of composites,
 and it enables us to use the technique of subspace substitution. This technique was introduced in \cite{Milton:2002:TOC}, Section 29.1, including the case 
 where  $\BS$ and $\BI-\BS$ project onto orthogonal subspaces and it was extended in \cite{Milton:2016:ETC}, Sections 7.8 and Chapter 8, to include the case where they project onto non-orthogonal subspaces.
 Define $\CU$, $\CE$, and $\CJ$ to be the three subspaces spanned by the three unit vectors 
 \beq \Bw_0=\bpm 1 \\0 \\ 0 \epm,\quad \Bw_1=\bpm 0 \\1 \\ 0 \epm,\quad \Bw_2=\bpm 0 \\0 \\ 1 \epm,
 \eeq{4.1A}
 respectively, so that $\BGG_i=\Bw_i\otimes\Bw_i$, $i=1,2,3$, are the projections onto $\CU$, $\CE$, and $\CJ$ respectively, i.e.
 \beq \BGG_0 = \bpm 1 & 0 & 0 \\ 0 & 0 &0 \\  0 & 0 & 0 \epm,\quad
 \BGG_1 = \bpm 0 & 0 & 0 \\ 0 & 1 &0 \\  0 & 0 & 0 \epm, \quad
 \BGG_2 = \bpm 0 & 0 & 0 \\ 0 & 0 &0 \\  0 & 0 & 1 \epm.
 \eeq{assproj}
 Then \eq{4.1} reduces to
 \beq J\Bw_0+ J_2\Bw_2=\widetilde{\BL}(E\Bw_0+ E_1\Bw_1)-h\Bw_0,
 \quad \widetilde{\BL}=\BI-\BS/\underline{z},
 \eeq{4.1B}
 which is a problem in the abstract theory of composites. More generally, 
 a problem in the abstract theory of composites
 takes the form: given $\BE_0\in\CU$, and a source term $\Bh$ in $\CH=\CU\oplus\CE\oplus\CJ$, and
 an operator $\widetilde{\BL}$ mapping $\CH$ to $\CH$, find $\BJ_0\in\CU$, $\BE\in\CE$ and $\BJ\in\CJ$ such that
 \beq \BJ_0+\BJ=\widetilde{\BL}(\BE_0+\BE)-\Bh. \eeq{etc}
 In our case, the subspaces $\CU$, $\CE$, and $\CJ$ are clearly orthogonal, but $\BS$ and $\BI-\BS$ do not generally project onto orthogonal subspaces
 when $s_1,s_2$ and $s_3$ are not all real.
 
 To find the norm of $\BS$ we consider its action on a possibly complex vector $\Ba$. We have
 \beq |\BS\Ba|=|\Bs(\Bs\cdot\Ba)|\leq |\Bs|^2|\Ba|, \eeq{4.1b}
 with equality when $\Ba$ is the complex conjugate of $\Bs$. Thus $\BS$ has norm 
 \beq \|\BS\|=|s_1|^2+|s_2|^2+|s_3|^2, \eeq{4.1c}
 and this will surely be greater than or equal to $1$ if \eq{4.1a} holds and $s_1,s_2$ and $s_3$ are either real or purely imaginary.
 For example, if $s_1$ is purely imaginary while $s_2$ and $s_3$ are purely real then \eq{4.1c} implies
 \beq 1=-|s_1|^2+|s_2|^2+|s_3|^2=\|\BS\|^2-2|s_1|^2, \eeq{4.1d}
 which forces $\|\BS\|$ to be greater than or equal to $1$. 
 
 The matrix equation \eq{4.1} is clearly satisfied with $E_1  =  (s_1s_2)E/(z-s_2^2)$ which gives the
 so called "effective equation"
 \beq  J =(1-1/z)E-h, \text{ where } z=\frac{\underline{z}-s_2^2}{s_1^2}.
 \eeq{4.4}
 Solving this last equation for $\underline{z}$ in terms of $z$ gives
 \beq \underline{z}=s_1^2z+s_2^2.
 \eeq{4.5a}

 Suppose now that in the extended abstract theory of composites
 we are interested in solving the equations
 \beq \BJ=\BL\BE-\Bh,\quad\text{with}\quad \BGG\BE=\BE,\quad\BGG\BJ=0,
 \eeq{4.6}
 where
 \beq \BL=\BI-\BQ/z, \eeq{4.6a}
 or equivalently in finding the resolvent \eq{1.1} with $\BA=\BGG\BQ\BGG$.
 Our preliminary linear algebra problem shows this is equivalent to solving
 \beq 
 \underbrace{\bpm
 	\BJ\\0 \\ \BJ_2 \epm}_{\underline{\BJ}} = \Biggl[\BI- \frac{1}{\underline{z}}
 \underbrace{\bpm s_1^2\BQ & s_1s_2\BQ & s_1s_3\BQ \\ s_1s_2\BQ &
 	s_2^2\BQ & s_2s_3\BQ \\  s_1s_3\BQ & s_2s_3\BQ & s_3^2\BQ \epm
 }_{\underline{\BQ}}\Biggr]\underbrace{\bpm \BE\\ \BE_1 \\ 0 \epm}_{\underline{\BE}}-
 \underbrace{\bpm \Bh\\ 0 \\ 0 \epm}_{\underline{\Bh}},
 \eeq{4.8}
 with $\underline{\BGG}\BE=\BE$ and $\underline{\BGG}\BJ=0$, in which
 \beq \underline{\BGG}=\bpm \BGG & 0 & 0 \cr 0 & \BI & 0 \cr 0 & 0 & 0 \epm.
 \eeq{4.8aa}
 Specifically, from \eq{4.8} we get 
 \beq
 \bpm
 \BQ\BJ\\0 \\ \BQ\BJ_2 \epm= \Biggl[\BI-\frac{1}{\underline{z}}
 \bpm s_1^2\BI & s_1s_2\BI & s_1s_3\BI \\ s_1s_2\BI &
 s_2^2\BI & s_2s_3\BI \\  s_1s_3\BI & s_2s_3\BI & s_3^2\BI \epm
 \Biggr]\bpm \BQ\BE\\ \BQ\BE_1 \\ 0 \epm-\bpm \BQ\Bh\\ 0 \\ 0 \epm,
 \eeq{4.8a}
 and 
 \beq 
 \bpm
(\BI-\BQ)\BJ\\0 \\ (\BI-\BQ)\BJ_2 \epm=\bpm (\BI-\BQ)\BE\\ (\BI-\BQ)\BE_1 \\ 0 \epm-\bpm (\BI-\BQ)\Bh\\ 0 \\ 0 \epm.
 \eeq{4.8b}
 Then, from \eq{4.8a}, \eq{4.1}, and \eq{4.4} we obtain
 \beq \BQ\BJ=(1-1/z)\BQ\BE-\BQ\Bh=\Gs\BQ\BE-\BQ\Bh, \eeq{4.8c}
 and \eq{4.8b} implies $(\BI-\BQ)\BJ=(\BI-\BQ)\BE-(\BI-\BQ)\Bh$. Together they imply 
 \eq{4.6} with $\BL$ given by \eq{4.6a}.

 We are back at an equivalent problem now involving a new resolvent. Specifically, since $\underline{\BJ}$ and
 $\underline{\BE}$ lie in orthogonal spaces, we have
 \beq \underline{\BJ}(\Bx)=\underline{\BL}(\Bx)\underline{\BE}(\Bx)
 -\underline{\Bh},\quad
 \underline{\BGG}_1\underline{\BE}=\underline{\BE},\quad \underline{\BGG}_1\underline{\BJ}=0,\quad
 \text{with}\quad \underline{\BL}=\BI-\underline{\BQ}/\underline{z}.
 \eeq{4.9}
 Let us now see how this can improve convergence. Applying $\underline{\BGG}$ to both sides of  the first equation in \eq{4.9} we get
 \beq \underline{\BE}=[\BI-\underline{\BGG}\underline{\BQ}/\underline{z}]^{-1}\underline{\Bh} =[\BI-(1-\underline{\Gs})\underline{\BGG}\underline{\BQ}]^{-1}\underline{\Bh}, 
 \eeq{4.11}
 in which
 \beq \underline{\BQ}=\bpm s_1^2\BQ & s_1s_2\BQ & s_1s_3\BQ \\ s_1s_2\BQ &
 s_2^2\BQ & s_2s_3\BQ \\  s_1s_3\BQ & s_2s_3\BQ & s_3^2\BQ \epm
 ,\quad
 \underline{\Gs}=1-1/\underline{z}=1+ \frac{1}{s_1^2z+s_2^2}. \eeq{4.12}
 Note that $\underline{\BQ}$ is a projection operator
 because both $\BS$ and $\BQ$ are projections and thus the inverse in \eq{4.11} has exactly the same form as in \eq{1.14A} with
 $z$ and $\BQ$ being replaced by $\underline{z}$ and  $\underline{\BQ}$ 
 Also $\underline{\sigma}=(\underline{z}-1)/\underline{z}$ can be re-expressed in the form
 \beq \underline{\sigma}=\frac{z-z^-}{z-z^+},
 \eeq{4.11a}
 with
 \beq z^+=-s_2^2/s_1^2, \quad z^-=-(1+s_2^2)/s_1^2.
 \eeq{4.11b}
 Note that $\underline{\sigma}=0$  and $\underline{\sigma}=\infty$ are obtained by substituting $z=z^-$ and  $z=z^+$ in \eq{4.12}. Given real $1>z^+>z^->0$ we need to
 choose $s_1$ and $s_2$ so that these equations are satisfied. This will
 necessitate complex solutions for $s_2$ since otherwise $z^+$
 will be negative. Explicitly, we have
 \beq s_1^2=\frac{1}{z^+ - z^-}, \quad  s_2^2=\frac{-z^+}{z^+ - z^-}, \eeq{4.12A}
 with $s_1$ being real and $s_2$ being purely imaginary,
 and so $\underline{\BQ}$  is no longer Hermitian implying $\underline{\BQ}$ and $\BI-\underline{\BQ}$ 
 project onto subspaces that are not orthogonal.
 This translates to a problem in the extended abstract theory of
 composites with a non-orthogonal subspace collection, as introduced in Chapter 8 of \cite{Milton:2016:ETC}. 
 
 Next, we follow the steps outlined in the previous section, though now the projection $\underline{\BQ}$ 
 does not have norm 1.  We obtain the expansion
 \beq [\BI-{\underline{\BGG}}{\underline{\BQ}}/\underline{z}]^{-1}= \underline{\BH}\sum_{n=0}^\infty v^n  \underline{\BD}^n
 \eeq{4.14A}
 in powers of
 \beq v=\frac{\sqrt{{\underline{\sigma}}}-1}
 {\sqrt{{\underline{\sigma}}}+1}=\frac{\sqrt{\frac{z-z^-}{z-z^+}}-1}
 {\sqrt{\frac{z-z^-}{z-z^+}}+1},
 \eeq{1.14AA}
 with
 \beq \underline{\BH}=2(\underline{\BL}+\sqrt{\underline{\Gs}}\BI)^{-1}=2\underline{\BQ/}(\underline{\Gs}+\sqrt{\underline{\Gs}})+2(\BI-\underline{\BQ})/(1+\sqrt{\underline{\Gs}}), \quad
   \underline{\BD}=(2\underline{\BQ}-\BI)(\BI-2\underline{\BGG}), \quad
 \eeq{4.16A}
  Note that the formula for $v$ is a composition of the maps
 \beq  v=\frac{w-1}{w+1},\quad w=\sqrt{\underline{\sigma}},\quad \underline{\sigma}=\frac{z-z^-}{z-z^+},
 \eeq{4.15}
 whose action is discussed in Section 3.
 The desired series expansion is obtained by substituting the above into the relation
 \beq [\BI-\BGG\BQ/z]^{-1}=
 \underline{\BG}^\dagger[\BI-\underline{\BGG}\underline{\BQ}/\underline{z}]^{-1} \underline{\BG},\quad
 \text{   with   } \underline{\BG}=\bpm \BI \cr 0 \cr 0 \epm,
 \eeq{4.16}
 giving
 \beq  [\BI-\BGG\BQ/z]^{-1}=\frac{2}{\sqrt{\underline{\sigma}}+1} \sum_{n=0}^\infty v^n\underline{\BG}^\dagger \underline{\BH}  \underline{\BD}^n\underline{\BG}.
 \eeq{4.16AAA}
 As this is the expansion of a self-adjoint operator in powers of $v$ we deduce that 
 the terms \newline $\underline{\BG}^\dagger \underline{\BH} \underline{\BD}^n\underline{\BG}$  must be self-adjoint 
 even though $ \underline{\BD}$ is not. 
 Since  $[\BI-\BGG\BQ/z]^{-1}$ has potentially singularities anywhere in the interval
 $[z^-,z^+]$, corresponding to $v$ being on the unit circle we deduce that the radius of convergence is at least $1$ so the rate of convergence is at least $\Gm_4=|v|$. To compare the  rates of convergence bounds $\Gm_2$ and $\Gm_4$ it is helpful to first express 
 both in terms of $w=\sqrt{\underline{\sigma}}$:
 \beq \Gm_2=\left|\frac{w^2-1}{w^2+1}\right|,
 \quad
 \Gm_4=\left|\frac{w-1}{w+1}\right|.
 \eeq{4.19}
 Thus the convergent rate of the first series with a given value $\underline{\sigma}=w$ will be the same as the second when  $\sigma=w^2$. 
 A similar analysis as that which lead to \eq{1.16aa} now implies
 \beq \Gm_2=\frac{2|v|}{|v^2+1|}, \quad \Gm_4=|v|. \eeq{4.20}
 Clearly if $v$ is in the unit disk then so is $v^2$ and $|v^2+1|$ is at most 2, equaling 2 when $v=\pm 1$. As the unit disk is the image of the complex
 $z$-plane, excluding the cut $z\in [z^-,z^+]$ we conclude that $\Gm_4\leq \Gm_2$.
 When $v$ is  small then $\Gm_4\approx 2\Gm_4$, both being very large. 
 We do not have a proof that $\Gm_4\leq \Gm_3$, although we suspect that to be the case.  The convergence has been 
 tested by Moulinec and Suquet (private communication) in a  model example of  a conducting composite and the convergence based on the series 
 \eq{4.16A} is better than that of the series \eq{1.14A}: see figure 8.8 in \cite{Milton:2016:ETC}. 
 
 Naturally, computing the action of $  \underline{\BD}$ on a field is more involved than 
 computing the action of $\BD$ on a field. So, if $z^-$ is close to zero and $z^+$ is close to 1, it is probably best to work with the series \eq{1.14A} rather than \eq{4.16A}. 
 
 \subsection{What if $z^+$ and $z^+$  are only known approximately?}
 
 We may not know $z^+$ and $z^+$ exactly but rather just have estimates 
  $z^+_e$ and $z^+_e$ for them. The treatment here follows that discussed in 
  Section 8.8 of \cite{Milton:2016:ETC}. Using the estimates gives
  \beq v=\left(
  \frac{\sqrt{\frac{z-z_e^-}{z-z_e^+}}-1}
  {\sqrt{\frac{z-z_e^-}{z-z_e^+}}+1}.
  \right)
 \eeq{F.1}
 The interesting case is when $z^+> z^+_e$ or $z^-_e>z^-$ since  otherwise 
 $z^-_e\BI\leq \BGG\BA\BGG \leq z^+_e\BI$ and all the previous results hold with 
 $z^+$ and $z^-$ replaced with $z^+_e$ and  $z^-_e$. Just looking at the case
 where $z^+> z^+_e$ and $z^-_e>z^-$ the points $z^-$ and $z^+$ map to
\beq v^-=\left(
\frac{\sqrt{\frac{z^--z_e^-}{z^--z_e^+}}-1}
{\sqrt{\frac{z^- -z_e^-}{z^- -z_e^+}}+1}
\right), \quad v^+=
\left(
\frac{\sqrt{\frac{z^+ -z_e^-}{z^+-z_e^+}}-1}
{\sqrt{\frac{z^+-z_e^-}{z^+-z_e^+}}+1}
\right).
 \eeq{F.2}
To the leading order in $z^--z_e^-$ and $z^+-z_e^+$ approximate to
 \beq
v^- \approx  \sqrt{\frac{z_e^--z^-}{z^+-z^-}}-1, \quad
	v^+\approx 1-
	\sqrt{\frac{z^+-z_e^+}{z^+ -z^-}}.
 \eeq{F.3}
 Assuming that the spectrum of $\BA$ includes both $z^-$ and $z^+$,
 the rate of convergence will be the greater of $|v/v^+|$ and $|v/v^-|$. Since \eq{F.3} implies
 $|v^- |\approx 1$ when $z^--z_e^-$ is small and  $|v^-+|\approx 1$ when
 $z^+-z_e^+$ is small  the rate of convergence will be close to that obtained when we know $z^+$ and $z^-$ exactly. 
 
 A natural way of obtain these estimates $z_e^-$ and $z_e^+$ is via the Rayleigh Ritz procedure. Thus, we pick appropriate finite dimensional subspaces $\CS^-\subset \BCE$ and $\CS^+\subset\BCE$ 
 and take
 \beqa z^-_e &= &\min_{\Ba\in \CS^- \,\, |\Ba|=1}(\BP^\dagger\BA\BP\Ba,\Ba), \nonum
 z^+_e & = &\max_{\Ba\in \CS^+,\,\,|\Ba|=1}(\BP^\dagger\BA\BP\Ba,\Ba),
 \eeqa{F.4}
 where $(\cdot,\cdot)$ denotes the inner product in $\BCE$ which we take to be the inner product in our Hilbert or vector space $\BCH$.  If, for example, the spectrum was discrete then ideally there should  be a field
 $\Ba\in \CS^-$ that is a good approximation to some field in the eigenspace
 of  $\BGG\BA\BGG$ having the lowest eigenvalue, and ideally there should be a field
 $\Ba\in \CS^+$ that is a good approximation to some field in the eigenspace
 of  $\BGG\BA\BGG$ having the maximum eigenvalue.
 
\section{Conclusion}

We have presented a novel approach to calculating 
a certain class of resolvents via series expansions. These series expansions provide a generally rapidly converging 
way of computing the action of these resolvents on fields
via an iterative method. While the method assumes the spectrum of $\BP^\dagger\BQ\BP$ lies inside a known,
or at least approximately known, interval $[z^-, z^+]$,
it seems likely that adaptive methods could be developed that do not require knowledge of $z^-$ and $z^+$.  The performance of the method improves upon some  methods,
and likely upon all the alternate methods presented here, at least in the asymptotic limit where one is close to convergence.
However, it remains to provide bounds on the error incurred when one truncates the series expansion. Clearly, explicit numerical investigations, beyond those of Moulinec and Suquet (private communication, 2015) as
summarized in figure 8.8  of \cite{Milton:2016:ETC}
 need to be made. These should compare the usefulness of the method with other available methods, both in terms of speed of convergence and computer memory required. In particular, it remains to compare the method presented here with more established methods, such as the conjugate gradient method (with $z$ real) or the biconjugate gradient method, biconjugate gradient stabilized method, or conjugate gradient squared method (with $z$ complex). When $z$ is real the condition number of the operator $(z\BI-\BP^\dagger\BQ\BP)$ we are inverting is
$\Gk=\underline{\Gs}=(z-z^-)/(z-z^+)$ when $z>z^+$ and
$\Gk=1/ \underline{\Gs}$ when $z<z^-$. So the convergence
rate of our new method 
\beq |v|=\left|\frac{\sqrt{\underline{\Gs}}-1}{\sqrt{\underline{\Gs}}+1}\right|
=\frac{\sqrt{\Gk}-1}{\sqrt{\Gk}+1}
\eeq{co1}
matches that of  the conjugate gradient method.
Comparisons of convergence between the  conjugate 
gradient method and the basic scheme \cite{Moulinec:1994:FNM} we discussed in section 3
have been made for conducting composites \cite{Zeman:2010:AFB} and elastic composites
\cite{Brisard:2010:FBM}. As far as memory goes the new method is likely to be more demanding than most other methods because $\underline{\BQ}$ and
$\underline{\BGG}$ act in spaces that have three times
the dimensionality of $\BCH$ when $\BCH$ is finite dimensional. 

Even if the conjugate gradient, biconjugate gradient, biconjugate gradient stabilized, and conjugate gradient squared methods perform better than the method ultimately presented here, the  idea of subspace substitution with non-orthogonal subspaces may have 
application in other areas of mathematics and numerical analysis in particular.

 \section*{Acknowledgments}
Matti Schneider is thanked for many helpful comments on the manuscript and 
Herve Moulinec and Pierre Suquet are thanked for additional comments.
The author is grateful to the National Science Foundation for support through Research Grant DMS-2107926.

\ifx \bblindex \undefined \def \bblindex #1{} \fi\ifx \bbljournal \undefined
\def \bbljournal #1{{\em #1}\index{#1@{\em #1}}} \fi\ifx \bblnumber
\undefined \def \bblnumber #1{{\bf #1}} \fi\ifx \bblvolume \undefined \def
\bblvolume #1{{\bf #1}} \fi\ifx \noopsort \undefined \def \noopsort #1{}
\fi\ifx \bblindex \undefined \def \bblindex #1{} \fi\ifx \bbljournal
\undefined \def \bbljournal #1{{\em #1}\index{#1@{\em #1}}} \fi\ifx
\bblnumber \undefined \def \bblnumber #1{{\bf #1}} \fi\ifx \bblvolume
\undefined \def \bblvolume #1{{\bf #1}} \fi\ifx \noopsort \undefined \def
\noopsort #1{} \fi

%\bibliographystyle{siam}
%\bibliography{/u/ma/milton/tcbook,/u/ma/milton/newref}
%\bibliography{/Users/milton/tcbook,/Users/milton/newref}

 \end{document}